\documentclass[runningheads,a4paper]{llncs}

\usepackage{amssymb}
\usepackage{amsmath}

\usepackage{graphicx}
\usepackage{comment}

\usepackage{algorithmicx}
\usepackage{algpseudocode}
\usepackage{algorithm}

\makeatletter
\newenvironment{breakablealgorithm}
  {
   \begin{center}
     \refstepcounter{algorithm}
     \hrule height.8pt depth0pt \kern2pt
     \renewcommand{\caption}[2][\relax]{
       {\raggedright\textbf{\fname@algorithm~\thealgorithm} ##2\par}%
       \ifx\relax##1\relax 
         \addcontentsline{loa}{algorithm}{\protect\numberline{\thealgorithm}##2}%
       \else 
         \addcontentsline{loa}{algorithm}{\protect\numberline{\thealgorithm}##1}%
       \fi
       \kern2pt\hrule\kern2pt
     }
  }{
     \kern2pt\hrule\relax
   \end{center}
  }
\makeatother

\usepackage{setspace}

\usepackage{listings}

\usepackage{url}
\usepackage{hyperref}

\usepackage{tabularx}
\usepackage{tabulary}
\usepackage{mathtools}

\usepackage{array,xcolor,colortbl}

\usepackage{multirow}
\usepackage{tensor}
\usepackage{hhline}

\newcolumntype{R}{>{\raggedleft}p{0.1\textwidth}}
\newcolumntype{C}{>{\centering\arraybackslash}p{0.2\textwidth}}
\def\XXX{\operatornamewithlimits{\mbox{\sf\Large X}}}

\usepackage{tikz}
\usetikzlibrary{external}
\urldef{\mailsa}\path|helmut.gfrerer@jku.at|
\urldef{\mailsb}\path|outrata@utia.cas.cz|
\urldef{\mailsc}\path|jvaldman@prf.jcu.cz |

\usepackage{hyperref}
\usepackage{mathtools}
\DeclarePairedDelimiter\ceil{\lceil}{\rceil}

\usepackage{amssymb}
\usepackage{tensor}

\newcommand{\R}{\mathbb{R}}
\newcommand{\tto}{\rightrightarrows}

\newcommand{\ssstar}{semismooth$^{*}$ }

\newcommand{\iter}[2]{\tensor*[^{#1}]{#2}{}}

\newcommand{\dom}{\mathrm{dom}\,}
\newcommand{\gph}{\mathrm{gph}\,}

\newcommand{\skalp}[1]{\langle #1\rangle}
\newcommand{\norm}[1]{\|#1\|}
\newcommand{\B}{{\cal B}}

\newcommand{\templ}[1]{\multicolumn{1}{|c}{#1}}
\newcommand{\tempr}[1]{\multicolumn{1}{c|}{#1}}
\newcommand{\templr}[1]{\multicolumn{1}{|c|}{#1}}

\begin{document}

\mainmatter  
\title{On the solution of contact problems with Tresca friction by the semismooth* Newton method}

\titlerunning{The semismooth* Newton method for friction contact problems}

\author{Helmut Gfrerer$^1$
and 
Ji\v r\'i V. Outrata$^2$ 
and 
Jan Valdman$^{2,3}$
\thanks{The work of the 2nd and the corresponding author was supported by the Czech Science Foundation (GACR), through the grant GF19-29646L.}
}

\authorrunning{Helmut Gfrerer and Ji\v r\'i V. Outrata and Jan Valdman}

\institute{$^1$ Institute of Computational Mathematics,  Johannes Kepler \\ University Linz, 
Linz, Austria;\\
Email: \mailsa\\ 
$^2$
Czech Academy of Sciences, Institute of Information Theory \\
and Automation, Prague, Czech Republic; \\
E-mail: \mailsb \\
$^3$
Institute of Mathematics, Faculty of Science, University \\
of South Bohemia, České Budějovice, Czech Republic; \\
Email: \mailsc
}
\maketitle
\begin{abstract}
An equilibrium of a linear elastic body subject to loading and satisfying the friction and contact conditions can be described by a variational inequality of the second kind and the respective discrete model attains the form of a generalized equation. To its numerical solution we apply the semismooth* Newton method by Gfrerer and Outrata (2019) in which, in contrast to most available Newton-type methods for inclusions, one approximates not only the single-valued but also the multi-valued part. This is performed on the basis of limiting (Morduchovich) coderivative. In our case of the Tresca friction, the multi-valued part amounts to the subdifferential of a convex function generated by the friction and contact conditions. The full 3D discrete problem is then reduced to the contact boundary. Implementation details of the semismooth* Newton method are provided and numerical tests demonstrate its superlinear
convergence and mesh independence.

\end{abstract}

\section{Introduction}

In \cite{GO} the authors developed a new, so-called \ssstar 
Newton-type method for the numerical solution of an inclusion 
\begin{equation}\label{1}
0\in F(x), 
\end{equation}
where $F:\R^n\tto\R^n$ is a closed-graph multifunction. In contrast to existing Newton-type method $F$ is approximated on the basis of the limiting (Mordukhovich) normal cone to the graph of $F$, computed at the respective point. Under appropriate assumptions, this method exhibits local superlinear convergence and, so far, it has been successfully implemented to the solution of a class of variational inequalities (VIs) of the first and second kind, cf. 
 \cite{GO} and \cite{GOV}. This contribution is devoted to the application of the \ssstar method to the discrete 3D contact problem with Tresca friction which is modelled as a VI of the second kind. Therefore the implementation can be conducted along the lines of \cite{GOV}. The paper has the following structure: In Section 2 we describe briefly the main conceptual iterative scheme of the method. Section 3 deals with the considered discrete contact problem and Section 4 concerns the suggested implementation. The results of numerical tests are then collected in Section 5. 
 
 We employ the following notation. For a cone $K$, $K^0$ stands for its (negative) polar and for a multifunction $F:\R^n\tto\R^n$, $\dom F$ and  
 $\gph F$ denote its domain and its graph, respectively. The symbol 
``$\xrightarrow{A}$'' means the convergence within the set $A$, $\norm{B}_F$ denotes the Frobenius norm of a matrix $B$ and $\B_\delta(x)$ signifies the $\delta-$ ball around $x$. 
 
 \section{The \ssstar Newton method}
 For the reader's convenience we recall fist the definition of the tangent cone and the limiting (Mordukhovich) normal cone. 
 \begin{definition}
 Let $A \subset \mathbb{R}^n$ be closed and $\bar x \in A$. Then 
 \begin{itemize}
     \item[(i)] the cone
     \[
T_{A}(\bar{x}):=\{u \in \mathbb{R}^{n} | \exists t_k \searrow 0, u_k \rightarrow u \mbox{ such that } \bar{x}+ t_ku_k \in A \, \forall k\}
\]
is called the {\em (Bouligand) tangent cone} to $A$ at $\bar x$;
\item[(ii)] the cone
\[ N_A(\bar x):= \{ x^* \in \mathbb{R}^n | \exists  x_k \xrightarrow{A} \bar{x}, x_k^* \rightarrow x^* \mbox{ such that } x_k^* \in (T_A(x_k))^0 \, \forall k  \}  \]
is called the {\em limiting (Mordukhovich) normal cone} to $A$ at $\bar x$.
 \end{itemize}
 \end{definition}

The latter cone will be extensively used in the sequel. Let us assign to a pair $(\tilde x, \tilde y) \in \gph F$ two $[n \times n]$ matrices $A, B$ such that their i-th rows, say $u_i^*, v_i^*$, fulfill the condition
\begin{equation}\label{2}
(u_i^*, -v_i^*) \in N_{\gph F}(\tilde x, \tilde y), \qquad i=1,2,\dots,n.  
\end{equation}
Moreover, let $\mathcal{A} F(\tilde x, \tilde y)$ be the set of matrices $A, B$ satisfying \eqref{2} and 
$$A_{reg} F(\tilde x, \tilde y) = \{ (A, B) \in \mathcal{A} F(\tilde x, \tilde y) |\, A \mbox{ is non-singular} \}.$$

The general conceptual iterative scheme of the \ssstar Newton method is stated in Algorithm \ref{AlgNewton} below.

\begin{breakablealgorithm}
\caption{Semismooth$^*$ Newton-type method for generalized equations}
\label{AlgNewton}
\hspace{0.2cm}
\begin{spacing}{1.2}
\begin{algorithmic}[1]
\State Choose a starting point $\iter{0}{x}$, set the iteration counter $k:=0$.
\State If $0\in F(\iter{k}{x})$, stop the algorithm.
\State Approximation step: compute 
\[(\hat x,\hat y)\in\gph F\] close to $(\iter{k}{x},0)$ such that ${\cal A}_{\rm reg}F(\hat x,\hat y)\not=\emptyset$.
\State Newton step: select
  $(A,B)\in {\cal A}_{\rm reg}F(\hat x,\hat y)$ and compute the new iterate
  \[ \iter{k+1}x=\hat x - A^{-1}B\hat y.\]
\State Set $k:=k+1$ and go to 2.
\end{algorithmic}
\end{spacing}
\end{breakablealgorithm}

Let $\bar x$ be a (local) solution of \eqref{1}. Since $\iter{k}{x}$ need not belong to $\dom F$ or $0$ need not be close to $F(\iter{k}{x})$ even if $\iter{k}{x}$ is close to $\bar x$; one performs in step 3 an approximate projection of $(\iter{k}{x}, 0)$ onto $\gph F$. Therefore the step 3 is called the 
{\em approximation step}. The {\em Newton step} 4 is related to the following fundamental property, according to which the method has been named.

\begin{definition}[\cite{GO}]
Let $(\tilde x, \tilde y) \in \gph F$. We say that $F$ is \ssstar at $(\tilde x, \tilde y)$ provided that for every $\epsilon>0$ there is some $\delta>0$ such that the inequality
\begin{equation}\label{4}
\vert 
\skalp{x^*,x- \tilde x} + \skalp{y^*,y- \tilde y}
\vert
\leq 
\epsilon \norm{(x,y)-(\tilde x,\tilde y)} \, \norm{(x^*,y^*)} 
\end{equation}
is valid for all $(x,y) \in \B_\delta(\tilde x, \tilde y)$ and for all $(x^*, y^*)
\in  N_{\gph F}(x,y).$
\end{definition}
If we assume that $F$ is \ssstar at $(\bar x, 0)$, then it follows from \eqref{4} that for every $\epsilon>0$ there is some $\delta>0$ such that for every $(x,y) \in \gph F \cap \B_\delta(\bar x, 0)  $ and every pair $(A,B) \in {\cal A}_{\rm reg}F(x,y)$ one has
$$ \norm{(x - A^{-1}B y) - \bar x} \leq \epsilon \norm{A^{-1}} \, \norm{(A \,\vdots\, B)}_F \, \norm{(x,y) - (\bar x, 0) }, $$
cf. \cite[Proposition 4.3]{GO}. This is the background for the Newton step in Algorithm \ref{AlgNewton}. 

Finally, concerning the convergence, assume that $F$ is \ssstar at $(\bar x, 0)$ and there are positive reals $L, \kappa$ such that for every $x \not \in F^{-1}(0)$ sufficiently close to $\bar x$ the set of quadruples $(\hat x, \hat y, A, B),$ satisfying the conditions 
\begin{eqnarray}
&&\norm{(\hat x - \bar x, \bar y)} \leq L \norm{ x-\bar x}, \label{5} \\
&&(A, B) \in {\cal A}_{\rm reg}F(\hat x,\hat y), \label{6} \\
&&\norm{A^{-1}} \, \norm{(A \,\vdots\, B)}_F \leq \kappa \label{7}
\end{eqnarray}
is nonempty. Then it follows from \cite[Theorem 4.4]{GO}, that Algorithm \ref{AlgNewton} either terminates at $\bar x$ after a finite number of steps or 
converges superlinearly to $\bar x$ whenever $\iter{0}{x}$ is sufficiently close to $\bar x$. 

The application of the \ssstar Newton methods to a concrete problem of type \eqref{1} requires thus the construction of an approximation step and the Newton step which fulfill conditions \eqref{5}-\eqref{7}. 

\section{The used model}

The fundamental results concerning unilateral contact problems with Coulomb friction have been established in \cite{NJH}. The infinite-dimensional model of the contact problem with Tresca friction in form of a variational inequality of the second kind can be found, e.g., in \cite{HHN,OKZ}. Other related friction-type contact problems are described, e.g., in \cite{NRV}.

\begin{figure}
\centering
 \begin{minipage}[c]{.45\textwidth}
\centering  
\includegraphics[width=0.99\textwidth]{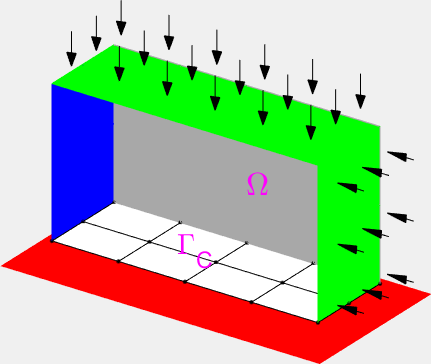}
\end{minipage}
\hspace{.06\textwidth}
\begin{minipage}[c]{.45\textwidth}
\vspace{0.68cm}
\centering  
\includegraphics[width=0.99\textwidth]{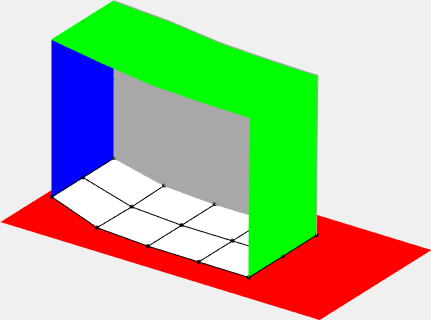}
\end{minipage}
\caption{The left picture depicts an undeformed elastic prism occupying domain $\Omega$ 
with the left (blue) face attached (Dirichlet condition) and some surface tractions applied to the right and top faces (depicted in green). They press the contact face $\Gamma_C$ against the (red) rigid plane foundation. Example of the resulting deformed body is depicted in the right picture. Front faces are not visualized. 
} \label{fig:setup}
\end{figure}

We assume that an elastic prism occupying domain $\Omega$ is pressed against a rigid plane foundation (cf. Figure \ref{fig:setup}). A full three-dimensional domain $\Omega$ is discretized by a mesh of brick elements and consists of $n$ nodes (vertices). The finite element method using trilinear basis functions is then applied to approximate a displacement field vector $u \in \mathbb{R}^{3 n}$ in each mesh node. Entries of $u$ are ordered in such a way that $u = (u^{1},u^{2},\ldots,u^{n})$ and 
the j-th node is associated with the pair $u^j=(u^{j}_{\tau}, u^{j}_{\nu}) \in \mathbb{R}^{2}\times \mathbb{R}$ of its {\em tangential} and {\em normal displacements}, respectively. 

A sparse stiffness matrix $K \in \mathbb{R}^{3 n \times 3 n}$ and the loading (column) vector $l \in \mathbb{R}^{3 n}$ are first assembled and then both condensed to incorporate zero displacements in Dirichlet nodes corresponding to the (blue) Dirichlet boundary. Secondly, all nodes not lying in the (bottom) contact face $\Gamma_{C}$ are eliminated by the {\em Schur complement} technique and the {\em Cholesky factorization} resulting in a dense matrix 
$\tilde{A} \in \mathbb{R}^{3 p \times 3 p}$ and a vector $\tilde{b} \in \mathbb{R}^{3 p}$, where $p \ll n$ is the number of $\Gamma_{C}$ nodes excluding Dirichlet boundary nodes. 

At last, all local $3 \times 3$ blocks of $\tilde A$ and all $3 \times 1$ blocks of $\tilde b$ are expanded to $4 \times 4$ blocks and $4 \times 1$ blocks, respectively
in order to incorporate the non-penetrability condition
\[ 0 \in u_{\nu}^i + N_{\mathbb{R}_{+}}(\lambda^{i}), \]
where $\lambda^i \in \mathbb{R}_{+}$ is the {\em Lagrange multiplier} associated  with non-penetrability constraint. Here and in the following, we assume that $i=1,\dots, p$. In this way, we obtain a dense regular matrix $A \in \mathbb{R}^{4 p \times 4 p}$ and a (column) vector $b \in \mathbb{R}^{4 p}$.

Finally, let us simplify the notation via
\[
\begin{split}
 x^{i}_{12}=(x^{i}_{1},x^{i}_{2})=u^{i}_{\tau} \in \mathbb{R}^{2}, \qquad x^{i}_{3} = u^{i}_{\nu} \in \mathbb{R}, \qquad x^{i}_{4} = \lambda^{i} \in \mathbb{R}
\end{split}
\]
to define a vector of unknowns  $x = (x^{1},x^{2},\ldots,x^{p}) \in \mathbb{R}^{4 p}$.

Following the development in \cite{BHKKO}, our model attains the form of generalized equation (GE) 
\begin{equation}\label{8}
0 \in f(x) + \widetilde{Q}(x),
\end{equation}
where the single-valued function 
$f: \mathbb{R}^{4p} \rightarrow \mathbb{R}^{4p}$ is given by 
\[ f(x)=A x - b\] and the multifunction
$\widetilde{Q}: \mathbb{R}^{4p} \rightrightarrows \mathbb{R}^{4p}$ by
\begin{equation*} 
\widetilde{Q}(x) = \XXX^{p}_{i=1} Q^i(x^{i}) \quad \mbox{  with  } Q^i(x^{i}) =
\left \{ \left [ \begin{array}{c}
-\phi \, \partial \| x^{i}_{12}\|\\
0 \\
N_{\mathbb{R}_{+}}(x^{i}_{4})
\end{array} 
\right ]   \right \}
\end{equation*}
with $\phi \geq 0$ being the {\em friction coefficient}.  
GEs of the type \eqref{8} have been studied in \cite{GOV} and so all theoretical results derived there are applicable. For our approach it is also important that the Jacobian $\nabla f(\bar x)$ is positive definite.

\section{Implementation of the \ssstar method}
In order to facilitate the approximation step we will solve, instead of GE \eqref{8}, the enhanced system

\begin{equation}\label{9}
0 \in \mathcal{F} (x,d) + \begin{bmatrix}
f(x) + \widetilde{Q}(d) \\
x - d
\end{bmatrix}
\end{equation}
in variables $(x,d) \in \mathbb{R}^{4p} \times \mathbb{R}^{4p}.$ Clearly $\bar x$ is a solution of \eqref{8}, if and only if $(\bar x, \bar x)$ is a solution of \eqref{9}.

In the approximation step we suggest to solve for all $i$ consecutively the next three low-dimensional strictly convex optimization problems:
 \begin{eqnarray*}
\mbox{(i)} &\underset{v \in \mathbb{R}^2}{\mbox{minimize}} &
\frac{1}{2}  \skalp{v,v} + \skalp{f^{i}_{12}(\iter{k}{x}),v} + \phi \norm{\iter{k}{x}_{12}^i + v}, \\
\mbox{(ii)} &\underset{v \in \mathbb{R}}{\mbox{minimize}} &
\frac{1}{2}  \skalp{v,v} + f^{i}_{3}(\iter{k}{x}) \cdot v, \\
\mbox{(iii)} &\underset{\iter{k}{x}_4^i + v \geq 0}{\mbox{minimize}} &
\frac{1}{2}  \skalp{v,v} + \iter{k}{x}_3^i \cdot v, 
\end{eqnarray*}
obtaining thus their unique solutions $\hat v^i_{12}, \hat v^i_{3}, \hat v^i_{4}$, respectively. The solutions can be ordered in vectors $\hat v^i = (\hat v^i_{12}, \hat v^i_{3}, \hat v^i_{4}) \in \mathbb{R}^4$ and all together in a vector 
\[ \hat v = (\hat v^1, \hat v^2, \dots,  \hat v^p) \in \mathbb{R}^{4p}. \]
Thereafter we compute the outcome of the approximation step via
$$\hat x = \iter{k}{x}, \qquad \hat d = \iter{k}{x} + \hat v, \qquad \hat y=(-\hat v, -\hat v).$$
Clearly $(\hat x, \hat d, \hat y) \in \gph \mathcal{F}$ and, using the theory \cite[Section 4]{GOV}, it is possible to show that condition \eqref{5} is fulfilled. 

In the Newton step we put
\begin{equation} \label{BChoice}
A=I, \qquad B = \begin{bmatrix} I  & & 0 \\ 0 \, & & G 
\end{bmatrix} D^{-1}, \qquad \mbox{ where    } D=\begin{bmatrix} \nabla f(\hat x) & -H \\ I & G
\end{bmatrix}.
\end{equation}
In \eqref{BChoice}, $I$ is an identity matrix and block diagonal matrices $G, H$ attain the form
\begin{equation}\label{12}
G=\mbox{diag}(G^1, G^2, \dots, G^p), \qquad H=\mbox{diag}(H^1, H^2, \dots, H^p),
\end{equation}
where the diagonal blocks $G^i, H^i \in \mathbb{R}^{4 \times 4}$ have the structure
\[
G^i=\left[\begin{array}{ccc}
\tempr{G^i_{1}} & \\ \cline{1-2}
 & \templr{1}&  \\  \cline{2-3}
 & & \templ{G^i_2}  
\end{array}\right], 
\qquad
H^i=\left[\begin{array}{ccc}
\tempr{H^i_1} & \\ \cline{1-2}
 & \templr{0}&  \\  \cline{2-3}
 & & \templ{H^i_2}  
\end{array}\right]
\]
and submatrices $G^i_1, H^i_1 \in \mathbb{R}^{2 \times 2}$ and scalar entries $G^i_2, H^i_2 \in \mathbb{R}$ are computed in dependence on values $\hat{d}^i_{12} \in \mathbb{R}^2$ and $\hat{d}^i_4 \in \mathbb{R}$ as follows:
\begin{itemize}
\item[$\bullet$] If $\hat{d}^i_{12}=0$ (sticking), we put $G^i_{1}=0, H^i_{1}=I$, otherwise we put  
\[ G^i_{1}=I, \qquad H^i_1 = \frac{\phi}{\norm{\hat{d}^i_{12}}^3}
\begin{bmatrix*}
(\hat{d}^i_2)^2 & -\hat{d}_1^i d_2^i \\  -\hat{d}_1^i \hat{d}_2^i & (\hat{d}^i_1)^2
\end{bmatrix*}.
\]
\item[$\bullet$] If $\hat{d}^i_4=0$ (no contact or weak contact), we put $G^i_{2}=0, H^i_{2}=1$, otherwise we put $G^i_{2}=1, H^i_{2}=0$.
\end{itemize}
This choice ensures that matrices $(I, B)$ with $B$ given by \eqref{BChoice} fulfill conditions \eqref{6},\eqref{7}  with $F$ replaced by $\mathcal{F}$. \\

\noindent {\bf Stopping rule} 
It is possible to show (even for more general Coulomb friction model \cite{BHKKO})  that there is a Lipschitz constant $c_L > 0$ such that,
$$\| (\hat{x},\hat{d})- (\bar{x},\bar{x}) \| \leq c_L \| \hat{y} \|,
$$
whenever the output of the approximation step lies in a sufficiently small neighborhood of $(\bar{u},\bar{u},0)$. It follows that, with a sufficiently small positive $\varepsilon$, the condition
\begin{equation}\label{stopping}
\| \hat{v}\| \leq \varepsilon,    
\end{equation}
tested after the approximation step, may serve as a simple yet efficient stopping rule. 

\subsubsection{Computational benchmark}

We assume that the domain $$\Omega = (0,2) \times (0,1) \times (0.1,1)$$ 
is described by elastic
parameters $E=2.1\cdot 10^9$ (Young's modulus), $\nu=0.277$ (Poisson's ratio) and subject to surface tractions
\begin{eqnarray}
&&f=(-5\cdot 10^8, 0, 0) \qquad \mbox{on the right-side face}, \\
&&f=(0, 0, -1\cdot 10^8 ) \qquad \mbox{on the top face}. 
\end{eqnarray}
and the friction coefficient $\phi=1$. 
\begin{table}[b]
\centering
\begin{tabular}{|c|rrr|rrr|}
\hline
level &  nodes &  assembly & Cholesky $\&$ & nodes & \multicolumn{2}{c|}{\ssstar solver}
\\
 $(\ell)$& 
$(n)$ & 
of K (sec)& Schur (sec)
& $(p)$ 
& time (sec) & iters \\
\hline
1 & 225 & 0.078 & 0.003  & 40& 0.017  & 5\\ 
2 & 225 & 0.031  & 0.003  & 40& 0.017  & 6\\ 
3 & 637 & 0.094 & 0.021  & 84& 0.047  & 6\\ 
4 & 1377 & 0.141 & 0.092  & 144& 0.101  & 6\\ 
5 & 4056 & 0.516 & 0.701  & 299& 0.507  & 7\\
6 & 9537 & 1.297 & 3.968  & 544& 1.928  & 7\\ 
7 & 27072 & 3.156 & 32.110  & 1104& 9.734  & 7\\
8 & 70785 & 18.672 & 1242.211  & 2112& 48.275  & 8\\
\hline
\end{tabular}
\caption{Performance of the MATLAB solver.} \label{performace}
\end{table}
The domain is uniformly divided into 
$e_x \cdot e_y \cdot e_z $ hexahedra (bricks), 
where 
$$e_x=\ceil{4 \cdot 2^{\ell/2}}, \quad e_y =\ceil{2 \cdot 2^{\ell/2}}, \quad e_z = \ceil{2 \cdot 2^{\ell/2}}$$ are numbers of hexahedra along with coordinate axis, $\ell$ denotes the mesh level of refinement and $\ceil{\cdot}$ the ceiling function. Consequently the number of $\Omega$ nodes $n$ and the number of $\Gamma_C$ nodes $p$ read
$$n(\ell)=(e_x+1) \cdot (e_y+1) \cdot (e_z+1), \qquad p(\ell)=e_x \cdot (e_y+1)$$ respectively.
\begin{figure}
\centering
 \begin{minipage}[c]{.46\textwidth}
\centering  
\includegraphics[width=0.99\textwidth]{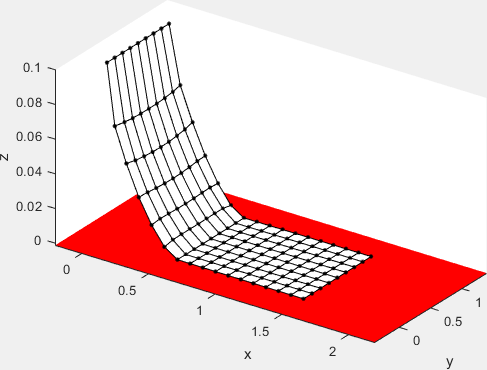}
\end{minipage}
\hspace{.06\textwidth}
\begin{minipage}[c]{.46\textwidth}
\centering  
\includegraphics[width=0.99\textwidth]{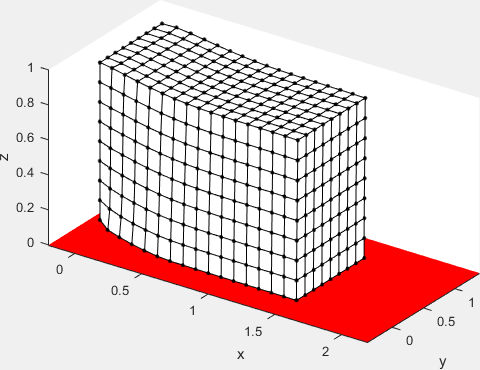}
\end{minipage}
\caption{The left picture depicts the deformed 
contact boundary and the right figure shows the corresponding deformed elastic prism, both pictures together with the (red) rigid plane foundation.
} \label{fig:solution}
\end{figure}
Table \ref{performace} reports on the performance of the whole method for various meshes assuming zero initial approximation $\iter{0}{x}=0$ and the stopping criterion 
$\epsilon=10^{-6}$. We can clearly see that the number of iterations of the \ssstar Newton method (displayed in the last column) only slightly increase with the mesh size. This behaviour shows that the method is mesh-independent. 

Figure \ref{fig:solution} visualizes a deformed contact boundary together with a deformation of the full domain $\Omega$ obtained by post-processing. Displacements of non-contact boundary nodes are then obtained from a linear system of equations with the matrix $K$ and the vector $l$. 

All pictures and running times were produced by our MATLAB code 
available for download and testing at
\begin{center}
\url{https://www.mathworks.com/matlabcentral/fileexchange/70255} .
\end{center}
It is based on original codes of \cite{BHKKO} and its performance is further enhanced by a vectorized assembly of $K$ using \cite{CSV}. 
\subsubsection{Concluding remarks and further perspectives}
The choice \eqref{BChoice} of matrices $A, B$ in the Newton step of the method is not unique and may be used to simplify the linear system in the Newton step. The convergence may be further accelerated by an appropriate scaling in the approximation step. 

\vspace{-0.3cm}
\subsubsection{Acknowledgment}
Authors are grateful to Petr Beremlijski (TU Ostrava) for providing original Matlab codes of \cite{BHKKO} and discussions leading to various improvements of our implementation. 

\vspace{-0.3cm}


\bibliographystyle{abbrv}

\end{document}